\documentclass[11pt,letterpaper]{amsart}
\usepackage[utf8x]{inputenc}

\usepackage{a4wide}
\usepackage{amssymb}
\usepackage{amsmath}
\usepackage{amsthm}
\usepackage[mathscr]{euscript}
\usepackage[colorlinks=true]{hyperref}
\usepackage[usenames,dvipsnames]{xcolor}
\usepackage{todonotes}
\usepackage{appendix}

\numberwithin{equation}{section}

\theoremstyle{plain}

\theoremstyle{definition}
\newtheorem{definition}{Definition}[section]

\newtheorem{Theorem}[definition]{Theorem}
\newtheorem{Proposition}[definition]{Proposition}
\newtheorem{Lemma}[definition]{Lemma}

\newtheorem{remark}[definition]{Remark}
\newtheorem{example}[definition]{Example}
\newtheorem{conjecture}[definition]{Conjecture}

\newcommand{\R}{\mathbb R}

\usepackage[xcolor]{changebar}
\cbcolor{blue}
\newcommand{\Ric}{\mathrm{Ric}}

\usepackage[inline]{enumitem}
\newcommand{\enumlabelformat}{\roman}
\newcommand{\enumlabelfont}[1]{#1}
\newlength{\thelabelsep}
\setlist{labelsep=\thelabelsep}
\setlist[enumerate,1]{font=\enumlabelfont,label=(\enumlabelformat*),leftmargin=2.5em}
\setlist[itemize]{leftmargin=2.5em,label=$-$}

\newcounter{inlineenum}
\renewcommand{\theinlineenum}{\enumlabelformat{inlineenum}}

\let\epsilon\varepsilon
\let\phi\varphi

\title{Cosmological rigidity via $p$-harmonic temporal functions}

\author{Argam Ohanyan}
\begin{document}

\date{\today}


\maketitle

\begin{abstract}

We prove a splitting theorem for cosmological spacetimes which admit a globally $p$-harmonic Cauchy temporal function for some $p < 1$, without making any use of timelike lines. In Bartnik's context of cosmological and timelike geodesically complete spacetimes, we establish a Liouville-type unboundedness result which shows that globally $p$-harmonic temporal functions, $p < 1$, are necessarily Cauchy temporal. Moreover, using the theory of elliptic PDEs, we show that locally Lipschitz $p$-harmonic temporal functions with $p < 1$ are in fact smooth. Altogether, our work reduces the resolution of Bartnik's cosmological splitting conjecture to the global weak solvability of the $p$-d'Alembert equation for some $p < 1$ in the class of locally Lipschitz temporal functions.

\vskip 1em

\medskip


\end{abstract}

\tableofcontents

\section{Introduction}
\label{section: introduction}

\subsection{Background}

The singularity theorems of Hawking and Penrose \cite{penrose1965gravitational, hawking1967occurrence, HawkingPenrose} show that singularities, i.e., incomplete causal geodesics, are a \emph{generic} occurrence in spacetimes. We refer to the standard texts \cite{HawkingEllis, ONeill, BeemEhrlichEasley} for their accessible treatment. Among the classical theorems, the Hawking--Penrose singularity theorem \cite{HawkingPenrose} is the most refined, makes the weakest a priori assumption on the global causality and is able to accommodate a variety of initial and boundary conditions to conclude the existence of singularities. In turn, it makes use of a global curvature nontriviality assumption called the \emph{genericity condition}, which requires nontrivial tidal forces to exist somewhere along every inextendible causal geodesic. In the special setting of \emph{cosmological} spacetimes, i.e., spacetimes with compact Cauchy surfaces satisfying the strong energy condition, the Hawking--Penrose theorem takes the following form:

\begin{Theorem}[Cosmological Hawking--Penrose singularity theorem]
Let $(M,g)$ be a globally hyperbolic spacetime with compact Cauchy surfaces satisfying the strong energy condition, i.e., $\Ric(v,v) \geq 0$ for all $v \in TM$ timelike, as well as the genericity condition, i.e., for every inextendible causal geodesic $\gamma$ there exists a parameter $t$ such that the tidal force operator $R:T_{\gamma(t)}M \to T_{\gamma(t)}M$, $v \mapsto R(v,\dot \gamma(t))\dot \gamma(t)$, is not the zero endomorphism. Then there exist incomplete causal geodesics in $M$.
\end{Theorem}

It is natural to ask just how generic the genericity condition actually is, i.e., what spacetimes remain complete in its absence? In 1988, Bartnik \cite{Bartnik1988remarks} famously conjectured that product spacetimes of the form $(\R \times \Sigma,dt^2 - h)$, where $(\Sigma,h)$ is a compact Riemannian manifold of nonnegative Ricci curvature, are the only complete cosmological spacetimes. More precisely:

\begin{conjecture}[Bartnik's cosmological splitting conjecture]
\label{conjecture: Bartnik}
Let $(M,g)$ be a globally hyperbolic spacetime with compact Cauchy surfaces satisfying the strong energy condition. If $(M,g)$ is timelike geodesically complete, then it is isometric to a product spacetime $(\R \times \Sigma,dt^2 - h)$, where $(\Sigma,h)$ is a compact Riemannian manifold of nonnegative Ricci curvature.
\end{conjecture}

The conjecture is known to hold under the stronger assumption of nonnegative timelike sectional curvature, and thus in spacetime dimension $2$ since timelike sectional and Ricci curvatures coincide, due to the work of Ehrlich--Galloway \cite{EhrlichGalloway}. Their approach has recently been extended to prove a Bartnik-type cosmological splitting theorem in the synthetic setting of Lorentzian length spaces \cite{FloresHerreraSolisBartnik}.

The principal difficulty is the following: Under the assumptions of Bartnik's conjecture, it is elementary to construct a causal line. If the latter is timelike, the claim follows from the Lorentzian splitting theorem \cite{eschenburg1988splitting, galloway1989lorentzian, newman1990proof}. Because of this observation, a variety of rather natural additional conditions have been imposed to guarantee timelikeness of the line, we refer to Galloway's excellent survey for a detailed discussion of these developments \cite{Gallowaysurvey}. While the plethora of contributions leads us to believe that Bartnik's splitting conjecture is correct, a number of interesting causality-theoretic properties have been recently observed for any potential counterexample to it by Minguzzi \cite{MinguzziWall}.

A different approach to Bartnik's conjecture is to find a CMC (i.e., constant mean curvature) Cauchy surface, which immediately implies the splitting \cite[Cor.\ 1]{Bartnik1988remarks}. However, in general, such CMC Cauchy surfaces are difficult to construct (cf.\ \cite{dilts2017spacetimes, GallowayLingSectionalCMC, GallowayLingRemarksCMC}), and there are examples (even in vacuum) of globally hyperbolic spacetimes with compact Cauchy surfaces which satisfy the strong energy condition but which do not contain any CMC Cauchy surfaces \cite{Bartnik1988remarks, ChruscielGallowayPollack, LingOhanyanExamples, OhanyanLingVacuum}. Let us mention the recent work of Galloway--Ling \cite{GallowayLingExpandingCMC} which shows the existence of a CMC Cauchy surface using mean curvature flow methods under the additional assumption that there exists an expanding Cauchy surface.

Aside from the CMC Cauchy surface approach to Bartnik's conjecture, far fewer works have tried to establish splitting directly without relying on the existence of a timelike line and thus the Lorentzian splitting theorem (for results of this flavor, see the works of Galloway--Vega \cite{GallowayVegaAchronal, GallowayVegaHausdorff}).

\subsection{Main results}

The aim of this article is to approach Bartnik's conjecture from an elliptic point of view, by making use of the degenerate elliptic $p$-d'Alembertian
\begin{equation}
    \square_p f:=-\mathrm{div}(|\nabla f|_g^{p-2} \nabla f),
\end{equation}
where $p <1$. This operator first appeared in the work of Mondino--Suhr \cite{MondinoSuhr} and its (degenerate) ellipticity properties were discovered and studied in a very general metric spacetime setting by Beran, Braun, Calisti, Gigli, McCann, Rott, Sämann and the author in \cite{octet} (see also Braun's followup work \cite{BraunExactdAlembertian}, the surveys by Braun \cite{braun2025perspectives} and McCann \cite{mccann2025survey} and the note by Galloway, McCann and the author \cite{GMO}). They were subsequently used by Braun, Gigli, McCann, Sämann and the author to give elliptic proofs of the Lorentzian splitting theorems \cite{Quintetsmooth} which are much closer to the proof of the Riemannian splitting theorem of Cheeger--Gromoll \cite{CheegerGromoll}. The strategy was adapted to low regularity metrics to prove a splitting theorem for ($C^1$-weighted) $C^1$-Lorentzian metrics by the same set of authors \cite{Quintetnonsmooth} and to Finsler spacetimes in order to establish a vastly improved splitting theorem by Caponio, Ohta and the author in \cite{CaponioOhanyanOhta}.

Let us outline our program and discuss our main results. First, we fix our terminology on temporal functions. We use the Lorentzian signature $(+,-,\dots,-)$ throughout, as indicated in Subsection \ref{subsection: notation and conventions}.

\begin{definition}[Temporal function]
Let $(M,g)$ be a spacetime, $X$ its time orientation vector field and $\Omega \subseteq M$ open. A locally Lipschitz function $f:\Omega \to \R$ is called \emph{temporal} if for every compact set $K \subseteq \Omega$ there exists $c_K > 0$ such that
\begin{equation}
    \min(g(\nabla f,\nabla f), g(X,\nabla f)) \geq c_K
\end{equation}
a.e.\ on $K$. If $\Omega = M$, a temporal function $f:M \to \R$ is called \emph{Cauchy temporal} if each level set $\{f = s\}$ is a Cauchy hypersurface.
\end{definition}

If $f \in C^1(\Omega)$, the definition reduces to the usual notion of temporal function which requires $\nabla f$ to be future timelike.

\begin{definition}[Weakly $p$-harmonic temporal function]
Let $(M,g)$ be a spacetime, $p < 1$, $\Omega \subseteq M$ open and $f:\Omega \to \R$ a locally Lipschitz temporal function. Then $f$ is called \emph{weakly $p$-harmonic} if
\begin{equation}
    \int_\Omega g(|\nabla f|_g^{p-2}\nabla f, \nabla \varphi) \, dvol_g = 0
\end{equation}
for every $\varphi \in C^\infty_c(\Omega)$.
\end{definition}

Our first main result, proved in Section \ref{section: smoothness}, establishes the smoothness of $p$-harmonic temporal functions:

\begin{Theorem}[Smoothness of $p$-harmonic temporal functions]
\label{Theorem: smoothness intro}
Let $(M,g)$ be a spacetime, $\Omega \subseteq M$ open, $p < 1$ and $f:\Omega \to \R$ a locally Lipschitz temporal function. If $f$ is weakly $p$-harmonic, then $f \in C^\infty(\Omega)$ and it is classically $p$-harmonic in the sense that
\begin{equation}
    \square_p f = -\mathrm{div}(|\nabla f|_g^{p-2} \nabla f) = 0.
\end{equation}
\end{Theorem}

Our arguments rely on results from the theory of elliptic PDEs \cite{GilbargTrudinger}: First, we use de Giorgi--Nash--Moser estimates to show that $f \in C^{1,\alpha}$ for some $\alpha \in (0,1)$, then Schauder estimates boost that regularity to $C^{2,\alpha}$. From here, classical elliptic regularity gives smoothness. When comparing to regularity results for solutions of the $p$-Laplace equation in the Riemannian setting \cite{Dibenedetto, lewisplaplacian}, where $C^{1,\alpha}$ can in general not be improved, we benefit from the fact that local boundedness and local uniform future timelikeness of $\nabla f$ restrict it to the domain of uniform ellipticity of $\square_p$, making a bootstrap all the way to smoothness possible. In the Riemannian setting, of course, $p$-harmonic functions are smooth on the open set where their gradient does not vanish, but the latter can in general not be excluded.

While solving the $p$-d'Alembert equation in not a part of this article, we refer to \cite{MinkowskipdAlemb}, where $p$-harmonic temporal functions constructed as minimizers of the Lorentzian $p$-energy on suitable domains in Minkowski spacetime are indeed shown to be temporal and locally Lipschitz, hence smooth by our Theorem \ref{Theorem: smoothness intro}.

Our next main result and certainly the one requiring the most substantial body of arguments is a Liouville-type unboundedness result for $p$-harmonic temporal functions in cosmological spacetimes which are timelike geodesically complete, i.e., precisely in Bartnik's setting.

\begin{Theorem}[Liouville-type theorem for $p$-harmonic temporal functions]
\label{Theorem: Liouville intro}
Let $(M,g)$ be a globally hyperbolic, timelike geodesically complete spacetime with compact Cauchy surfaces satisfying the strong energy condition. If $f: M \to \R$ is a $p$-harmonic smooth temporal function, then $f$ is Cauchy temporal.
\end{Theorem}

The proof, which is the content of Section \ref{section: Liouville}, is lengthy and involved. It makes full use of every assumption made (see Subsection \ref{subsection: necessity of assumptions}). A useful tool in the proof, also relevant for the next main result, is the divergence theorem, which, due to compactness of Cauchy surfaces, allows us to fix a background Cauchy temporal function and exploit divergence-freeness of the vector field $|\nabla|_g^{p-2} \nabla f$ in order to reduce spacetime integrals to integrals on level sets of that function. We refer to \cite[Thm.\ 1.3.16]{BarLorentzianwave} for a version of the Lorentzian divergence theorem that we will use in various places.

Our final main result, established in Section \ref{section: cosmological rigidity}, is a cosmological splitting theorem under the assumption that there exists a globally $p$-harmonic Cauchy temporal function.

\begin{Theorem}[Cosmological splitting via $p$-harmonicity]
\label{Theorem: cosmological splitting intro}
Let $(M,g)$ be a globally hyperbolic spacetime with compact Cauchy surfaces satisfying the strong energy condition. If there exists a $p$-harmonic Cauchy temporal function $f:M \to \R$, then $(M,g)$ splits isometrically as $(\R \times \Sigma, dt^2 - h)$, where $(\Sigma,h)$ is a compact Riemannian manifold of nonnegative Ricci curvature. In fact, $\Sigma = \{f = 0\}$ and the splitting map is the flow of $\nabla f$ (upon normalization).
\end{Theorem}

In the elliptic proof of the Lorentzian splitting theorems \cite{Quintetsmooth}, a $p$-harmonic temporal function is provided by the Busemann function associated with the given complete timelike line, whose gradient is of unit norm. Consequently, the $p$-Bochner--Ohta identity (see Theorem \ref{Theorem: p Bochner Ohta} below) immediately implied that the Hessian of the Busemann function vanishes. In our cosmological splitting result above, no assumption on $|\nabla f|_g$ is necessary. The aforementioned divergence-freeness of the vector field $J$ together with compactness of Cauchy surfaces allows us to show global concavity of a certain positive function constructed from integrals along the level sets of $f$, which are thus constant and this leads to the vanishing of the additional term in the $p$-Bochner--Ohta identity, allowing us to conclude that $\mathrm{Hess}\, f = 0$. From here, it is easy to conclude that the flow of $\nabla f$ (upon normalizing by the a posteriori constant $|\nabla f|_g$) must be the desired isometry to a product spacetime.

Altogether, our main results allow us to give the following equivalent reformulation of Bartnik's conjecture in terms of the existence of global weak solutions of the $p$-d'Alembert equation:

\begin{Theorem}[Bartnik's conjecture $\Leftrightarrow$ weak solvability of $p$-d'Alembert equation]
\label{Theorem: Bartnik equivalent intro}
Let $(M,g)$ be a cosmological, timelike geodesically complete spacetime with compact Cauchy surfaces satisfying the strong energy condition. Then the following are equivalent:
\begin{enumerate}
    \item $(M,g)$ is isometric to a product spacetime $(\R \times \Sigma,dt^2-h)$, where $(\Sigma,h)$ is a compact Riemannian manifold with nonnegative Ricci curvature.
    \item There exists $p < 1$ and a locally Lipschitz temporal function $f:M\to \R$ that is weakly $p$-harmonic.
\end{enumerate}
\end{Theorem}

The direction (ii)$\Rightarrow$(i) follows from our main results, i.e., Theorems \ref{Theorem: smoothness intro}, \ref{Theorem: Liouville intro} and \ref{Theorem: cosmological splitting intro}. The direction (i)$\Rightarrow$(ii) is of course trivial, since the projection onto $\R$ is $p$-harmonic for every $p$ in a product spacetime.

\subsection{Discussion and outlook}

Based on Theorem \ref{Theorem: Bartnik equivalent intro}, let us formulate the following conjecture, which is equivalent to Bartnik's:

\begin{conjecture}
\label{conjecture: main conjecture p harmonic}
Let $(M,g)$ be a globally hyperbolic, timelike geodesically complete spacetime with compact Cauchy surfaces satisfying the strong energy condition. Then there exists $p < 1$ and a locally Lipschitz temporal function $f:M\to \R$ which is weakly $p$-harmonic.
\end{conjecture}

Conjecture \ref{conjecture: main conjecture p harmonic} may prove to be a simpler approach to Bartnik's conjecture in comparison with previous ones. We remark that naive Busemann-like approaches to construct $p$-harmonic functions do not work in Bartnik's setting, see the work of Galloway, McCann and the author \cite{GMO} for a detailed discussion. On suitable local domains, such as between level sets of a fixed background Cauchy temporal function, we expect to be able to find $p$-harmonic temporal functions as minimizers of the Lorentzian $p$-energy $\tilde f \mapsto \int_\Omega H_p(d\tilde f) \, dvol_g$, where $H_p$ is the $p$-Hamiltonian (cf. Definition \ref{def: p-hamiltonian}) subject to certain well-chosen boundary conditions, see \cite{MinkowskipdAlemb} for results in this direction on Minkowski spacetime. The major remaining challenge is to construct all of these local $p$-harmonic temporal functions in a way that their temporality and local Lipschitzness is uniform, as one can then pass to the limit as follows:

\begin{Proposition}[Global $p$-harmonicity in the limit]
\label{Proposition: global p harmonicity in the limit}
Let $(M,g)$ be a spacetime, $p  < 1$, $(\Omega_n)$ a nested exhaustion of $M$ by open sets, and $f_n:\Omega_n \to \R$ $p$-harmonic. If $(f_n)$ is locally uniformly Lipschitz and $(\nabla f_n)$ is locally uniformly future timelike, i.e., for every $K \subseteq M$ compact there exists $c_K > 0$ such that
\begin{equation}
    \inf_{n \geq n_K}\min(g(\nabla f_n,\nabla f_n), g(X,\nabla f_n)) \geq c_K,
\end{equation}
where $n_K$ is such that $K \subseteq \Omega_{n_K}$ and $X$ is the time orientation, then for every $o \in M$ there exists a $p$-harmonic temporal function $f:M \to \R$ such that, up to subsequences, $f_n - f_n(o)\to f$ in $C^1_{loc}$.
\end{Proposition}
\begin{proof}
    We can fix some $o \in M$ (w.l.o.g.\ $o \in \Omega_1$) and assume that $f_n(o) = 0$, otherwise replace $f_n$ by $f_n - f_n(o)$, which are still locally uniformly Lipschitz and locally uniformly temporal functions that are $p$-harmonic on $\Omega_n$.

    The arguments given in the proof of Theorem \ref{Theorem: smoothness}, specifically \eqref{eq: alpha Holder estimate on diff quot}, allow us to show that $(df_n)$ is locally uniformly $\alpha$-Hölder for some $\alpha \in (0,1)$ as a consequence of local uniform Lipschitzness, local uniform temporality and $p$-harmonicity. Thus, we may apply the Arzéla--Ascoli theorem to conclude the existence of a $C^1$-function $f: M \to \R$ such that, up to subsequences $f_n \to f$ in $C^1_{loc}$. Local uniform future timelikeness of $(\nabla f_n)$ implies that $\nabla f$ is future timelike and bounded, i.e., $f$ is temporal. The weak formulation of $p$-harmonicity is easily seen to pass to the limit under $C^1_{loc}$-convergence, so we conclude that $f$ is $p$-harmonic.
\end{proof}

Let us formulate the existence of the sequence $(f_n)$ as above explicitly as a conjecture, which is what we think will be the result to be targeted if one approaches Bartnik's conjecture from the point of view we have presented here.

\begin{conjecture}
\label{conjecture: locally uniform p harmonic sequence}
    Let $(M,g)$ be a globally hyperbolic, timelike geodesically complete spacetime with compact Cauchy surfaces satisfying the strong energy condition. Then there exists $p < 1$, an exhaustion $(\Omega_n)$ of $M$ by nested open sets, as well as $f_n:\Omega_n \to \R$ $p$-harmonic such that $(f_n)$ is locally uniformly Lipschitz and $(\nabla f_n)$ is locally uniformly future timelike.
\end{conjecture}

We expect a reasonable choice of $\Omega_n$ to be $I^+(\Sigma_{-n}) \cap I^-(\Sigma_n)$, where $\Sigma_s = \{\tau = s\}$ are the level sets of some background Cauchy temporal function. As a consequence of Theorem \ref{Theorem: Bartnik equivalent intro} and Proposition \ref{Proposition: global p harmonicity in the limit}, Conjectures \ref{conjecture: Bartnik}, \ref{conjecture: main conjecture p harmonic} and \ref{conjecture: locally uniform p harmonic sequence} are all equivalent to one another.

\begin{remark}[Extension to $u$-harmonic functions]
As a consequence of the work of Sálamo Candal and the author \cite{ohanyansalamocandal}, the $p$-d'Alembertian operators for $p < 1$ can be viewed as a subclass of a much larger class of degenerate elliptic operators, namely the $u$-d'Alembertian operators
\begin{equation}
    \square_u f:=-\mathrm{div}\left(\frac{u'(|\nabla f|_g)}{|\nabla f|_g} \nabla f\right),
\end{equation}
where $u:(0,\infty) \to \R$ is such that $u \in C^\infty$, $u' > 0$, $u'' < 0$, and $u':(0,\infty) \to (0,\infty)$ is surjective. For $u = u_p$, $p < 1$, this definition recovers the standard $p$-d'Alembertian. It would be interesting to investigate whether our results extend to $u$-harmonic temporal functions. Given that our arguments only rely on the properties of $u_p$ shared by this more general class of functions $u$, and not on the $p$-homogeneity of $u_p$, we expect this to be the case.
\end{remark}

\subsection{Notation and conventions}
\label{subsection: notation and conventions}

Smooth manifolds are by convention connected, 2nd countable and Hausdorff. A Lorentzian metric on a smooth manifold is smooth and has signature $(+,-,\dots,-)$. The symbol $g^*$ will denote the \emph{inverse metric} of $g$ acting on covectors. A vector $v$ is timelike if $g(v,v) > 0$, causal if $g(v,v) \geq 0$. Timelikeness and causality of covectors is defined with respect to $g^*$ in the analogous way. We write $|v|_g:=\sqrt{g(v,v)}$ for the Lorentzian norm of a causal vector $v$, $|\alpha|_{g^*}$ for the Lorentzian norm of a causal covector $\alpha$, and $\mathrm{Hess} \, f$ for the $(0,2)$ $g$-Hessian of a function $f$. Given a covector $\alpha$, we denote by $\alpha^{\sharp} := g^*(\alpha,\cdot)$ the corresponding vector assigned to it via the $g$-musical isomorphisms. Our convention for curvature is $R(X,Y)Z = \nabla_X \nabla_Y Z - \nabla_Y \nabla_X Z - \nabla_{[X,Y]}Z$ for vector fields $X,Y,Z$. By $vol_g$ we denote the volume measure induced by $g$ (which is not necessarily a volume form, since we do not assume orientability). Given a real-valued function $f$, we write $f_+:=\max(f,0)$ for its positive part.

\section{Preliminaries}

We recall some preliminary notions and results concerning the $p$-Hamiltonian \cite{MinguzziConestructures, Mccanndisplacement} and the $p$-d'Alembertian \cite{octet, Quintetsmooth}.

\begin{definition}[Cauchy hypersurface]
Let $(M,g)$ be a spacetime. A subset $\Sigma \subseteq M$ is called a \emph{Cauchy (hyper)surface} if it is met uniquely by every inextendible causal curve.
\end{definition}

\begin{definition}[(Cauchy) temporal function]
Let $(M,g)$ be a spacetime, $X \in \mathfrak X(M)$ the time orientation vector field, and $\Omega \subseteq M$ open. A locally Lipschitz function $f:\Omega \to \R$ is called \emph{temporal} if for every compact set $K \subseteq \Omega$ there exists $c_K > 0$ such that a.e.\ on $K$,
\begin{equation}
    \min(g(\nabla f, \nabla f), g(X,\nabla f)) \geq c_K.
\end{equation}
If $\Omega = M$, a temporal function $f:M \to \R$ is called \emph{Cauchy temporal} if each level set $\{f = s\}$, $s \in \R$, is a Cauchy hypersurface.
\end{definition}

If $f$ is a Cauchy temporal function of regularity $C^k$, $k  \geq 1$, then its level sets are spacelike Cauchy hypersurfaces of regularity $C^k$. Let us recall the following useful reformulation of Cauchy temporality: A temporal function is Cauchy temporal if and only it is unbounded along inextendible causal curves. Precisely, this means that whenever $\gamma:(a,b) \to M$ is a future and past inextendible causal curve, then $(f \circ \gamma)(t) \to +\infty$ if $t \uparrow b$ and $(f \circ \gamma)(t) \to - \infty$ if $t \downarrow a$.

\begin{definition}[$p$-harmonic]
Let $(M,g)$ be a spacetime, $p < 1$ and $\Omega \subseteq M$ open. A locally Lipschitz temporal function $f:\Omega \to \R$ is called \emph{$p$-harmonic} if for every $\varphi \in C^\infty_c(\Omega)$,
\begin{equation}
    \int_\Omega g(|\nabla f|_g^{p-2} \nabla f, \nabla \varphi) \, \mathrm{dvol}_g = 0.
\end{equation}
\end{definition}
\begin{remark}[Classical $p$-harmonicity]
If $f$ is $C^2$ and temporal, then $p$-harmonicity is equivalent to $\square_p f = 0$, where
\begin{equation}
    \square_p f:=-\mathrm{div}(|\nabla f|_g^{p-2} \nabla f)
\end{equation}
is the \emph{$p$-d'Alembertian}.
\end{remark}
\begin{definition}[$p$-Hamiltonian]
\label{def: p-hamiltonian}
Given $p < 1$, we define the \emph{$p$-Hamiltonian} $H_p:T^*M \to \R \cup \{+\infty\}$ via
\begin{equation}
    H_p(\omega):= \begin{cases}- u_p(|\omega|_{g^*}), & \omega \text{ future timelike},\\ +\infty & \text{otherwise}. \end{cases}
\end{equation}
Here $u_p(t):=p^{-1} t^p$ if $p \neq 0$ and $u_0(x):=\log(x)$.
\end{definition}
\begin{remark}[Fiberwise derivatives of $H_p$]
For later reference, we recall the first and second fiberwise derivatives of $H_p$, $p < 1$:
\begin{align}
    D H_p|_\omega \equiv DH_p(\omega) &= -|\omega|_{g^*}^{p-2} \omega^\sharp,\\
    D^2H_p|_{\omega} \equiv D^2H_p(\omega) &= |\omega|_{g^*}^{p-2}\left((2-p) \frac{\omega^\sharp  \otimes \omega^\sharp}{|\omega|_{g^*}^2} - g^*\right). \label{eq: D^2 H}
\end{align}
We emphasize that for each covector $\omega$, the application of $D^2H_p(\omega)$ to a covector $\sigma$ produces a vector.
\end{remark}
The next result is the central ingredient in the proof of our cosmological splitting theorem, see \cite[Lem.\ 12]{Quintetsmooth}.

\begin{Theorem}[$p$-Bochner--Ohta identity]
\label{Theorem: p Bochner Ohta}
Let $(M,g)$ be a spacetime, $f \in C^3(M)$ temporal and $p < 1$. Then
\begin{equation}
    -\mathrm{div}\left(D^2H_p|_{df} \bigg[d(u_p(|df|_{g^*}))\bigg]\right) - DH_p|_{df} [d \square_p f] = |\mathrm{Hess} \, f|_{D^2H_p|_{df}}^2 + \mathrm{Ric}(DH_p|_{df},DH_p|_{df}),
\end{equation}
where 
\begin{equation}
    |\mathrm{Hess} \,f|^2_{D^2 H_p|_{df}} :=\mathrm{Tr}(D^2H_p|_{df} (\mathrm{Hess} \,  f) D^2H_p|_{df} (\mathrm{Hess} \, f)) \geq 0
\end{equation}
and it is equal to $0$ if and only if $\mathrm{Hess} \, f = 0$.
\end{Theorem}

\section{Smoothness of $p$-harmonic temporal functions}
\label{section: smoothness}

The aim of this section is to show that $p$-harmonic temporal functions are smooth. More precisely:

\begin{Theorem}[$p$-harmonicity $\Rightarrow$ smoothness]
\label{Theorem: smoothness}
Let $(M,g)$ be a spacetime, $p<1$, $\Omega \subseteq M$ open. If $f:\Omega \to \R$ is a locally Lipschitz temporal function (i.e., $\nabla f$ is locally uniformly future timelike a.e.\ on $\Omega$) that is weakly $p$-harmonic, then $f$ is smooth and is classically $p$-harmonic, i.e., on $\Omega$,
\begin{equation}
    \square_p f = - \mathrm{div}(|\nabla f|_g^{p-2} \nabla f) = 0.
\end{equation}
\end{Theorem}
\begin{proof}
    This is a local statement, so we may assume that $\Omega = M = \R^n$, $(\R^n,g)$ is a globally hyperbolic spacetime, and that $f$ is a global weak solution of $\square_p f = 0$. We will use the notation
    \begin{equation}
        Q(x;\xi):=g^{ij}(x)\xi_i\xi_j, \quad \mu(x):=\sqrt{|\det g(x)|}.
    \end{equation}
    Then the $p$-d'Alembert equation for some function $\tilde f$ may be rewritten as
    \begin{equation}
        \partial_i F^i(x;D\tilde f) = 0, \quad \text{where} \quad F^i(x,\xi):=-\mu(x) Q(x;\xi)^{(p-2)/2}g^{ij}(x)\xi_j.
    \end{equation}
    If $x \in \R^n$ and $\xi \in T_x^*\R^n$ is a timelike covector, letting $w:=\sqrt{Q(x;\xi)}$ and $N^i:=\xi^i/w$, where $\xi^i:=g^{ij}\xi_j$, we have that
    \begin{equation}
        \frac{\partial F^i}{\partial \xi_j}(x;\xi) = \mu w^{p-2}((2-p)N^i N^j - g^{ij}).
    \end{equation}
    This is easily seen to be positive definite: If $\eta \neq 0$ is a covector, then
    \begin{equation}
        \frac{\partial F^i}{\partial \xi_j}(x;\xi) \eta_i \eta_j > 0.
    \end{equation}
    Moreover, for any compact subset $\mathcal K$ of the future timelike covectors in $T^* \R^n$ there exist $0 < \lambda \leq \Lambda$ such that for all $(x,\xi) \in \mathcal K$ and all covectors $\eta$,
    \begin{equation}
        \lambda |\eta|^2 \leq \frac{\partial F^i}{\partial \xi_j}(x;\xi) \eta_i \eta_j \leq \Lambda |\eta|^2,
    \end{equation}
    where $|\eta|$ denotes the Euclidean norm of $\eta$. Choose a sufficiently small coordinate ball $B_{2r}$ (w.l.o.g.\ around $0$). By local Lipschitzness and uniform temporality, after possibly
shrinking the ball, there exist $L,\delta>0$ such that
\begin{equation}
\label{eq: Df Q}
|Df(x)|\le L,\qquad
Q(z;Df(x))\ge\delta,\qquad
Df(x)[X(z)]\ge\delta
\end{equation}
for almost every $x\in B_{2r}$ and every $z\in\overline B_{2r}$, where $X$ is the time orientation vector field. The triangle and reverse triangle inequalities preserves these bounds under convex
combinations of the coordinate covectors. Thus, with
$\xi_t=(1-t)Df(x)+tDf(y)$, these covectors remain in a common compact
subset of the set of future timelike covectors, and there exist
$0<\lambda\le\Lambda$ such that
\begin{equation}
\lambda|\eta|^2
\le \frac{\partial F^i}{\partial\xi_j}(z;\xi_t)\eta_i\eta_j
\le \Lambda|\eta|^2
\end{equation}
uniformly for all $x,y\in E$, $z\in\overline{B}_{2r}$,
$t\in[0,1]$ and $\eta\in\mathbb R^n$, where
$E\subseteq B_{2r}$ is a fixed set of full measure on which
$Df$ is defined and \eqref{eq: Df Q} holds.

    For $0 < |h| < r$ consider the finite difference quotients
    \begin{equation}
        v_h(x):=\frac{f(x+he_k) - f(x)}{h},
    \end{equation}
    where $(e_1,\dots,e_n)$ is the standard basis of $\R^n$. Each $v_h$ is locally Lipschitz with a.e.\ derivative
    \begin{equation}
        \partial_j v_h(x) = \frac{\partial_jf(x+ he_k) - \partial_j f(x)}{h}.
    \end{equation}
    By weak $p$-harmonicity of $f$, for every $\varphi \in C^\infty_c(\R^n)$ we have
    \begin{equation}
        \int_{\R^n} F^i(x;Df(x)) \, \partial_i\varphi(x) \, d x = 0.
    \end{equation}
    Testing weak $p$-harmonicity against $\tilde \varphi(x):=\varphi(x - he_k)$, changing variables and subtracting the above yields
    \begin{equation}
    \label{eq: diff quot weak identity}
        \int_{B_r} \frac{F^i(x+ he_k; Df(x+he_k)) -F^i(x;Df(x) )}{h} \, \partial_i \varphi(x) \, dx = 0
    \end{equation}
    for every $\varphi \in C^\infty_c(B_r)$. We may write
    \begin{align}
        &F^i(x+he_k;Df(x+he_k)) - F^i(x;Df(x)) = \\&F^i(x+he_k;Df(x+he_k)) - F^i(x+he_k; Df(x))  \\&+ \underbrace{F^i(x+he_k; Df(x))- F^i(x;Df(x))}_{=:hb^i_h(x)}
    \end{align}
    The first term on the right hand side, upon division by $h$, can be rewritten using the fundamental theorem of calculus as $a^{ij}_h(x) \partial_j v_h(x)$, where
    \begin{equation}
        a^{ij}_h(x):=\int_0^1 \frac{\partial}{\partial \xi_j} F^i(x+he_k;(1-t)Df(x) + tDf(x+he_k)) \, dt.
    \end{equation}
    This is similar to how the $p$-d'Alembertian was linearized to apply the maximum principle to the Busemann functions in the proof of \cite[Prop.\ 9]{Quintetsmooth}. Substituting into \eqref{eq: diff quot weak identity} yields
    \begin{equation}
        \int_{B_r} (a^{ij}_h(x)\partial_jv_h(x) + b^i_h(x)) \, \partial_i\varphi(x) = 0
    \end{equation}
    for every $\varphi \in C^\infty_c(B_r)$, i.e., we have that
    \begin{equation}
    \label{eq: alpha Holder estimate on diff quot}
        \partial_i (a^{ij}_h(x)\partial_j v_h(x) + b^i_h(x)) = 0
    \end{equation}
    weakly on $B_r$. The uniform ellipticity of $a_h^{ij}$ (uniform also in $h$) was already argued in \cite{Quintetsmooth}, and it is easily seen (as a consequence of the fact that $Df$ stays in a compact subset of the future timelike cone in $T^*B_r$) that $a_h^{ij}$ and $b^i_h$ are bounded on $B_r$ independently of $h$. Since $v_h \in W^{1,\infty}(B_r) \subseteq W^{1,2}(B_r)$, we can apply the inhomogeneous De Giorgi--Nash--Moser estimate \cite[Thm.\ 8.24]{GilbargTrudinger} to see that there exists an $\alpha \in (0,1)$ such that
    \begin{equation}
        \|v_h\|_{C^{0,\alpha}(B_{r/2})} \leq C_*,
    \end{equation}
    where $C_* > 0$ is independent of $h$ (and $k$). Passing to the limit $h \to 0$ allows us to conclude that $f \in C^{1,\alpha}(B_{r/2})$.

    Next, we will show that $f$ is in fact $C^{2,\alpha}$ on some smaller ball. Let $0 < |h| < r/8$. Given $x,y \in B_{r/4}$, we have (using the fact that all of the covectors considered stay inside a compact subset of the future timelike covectors) that
    \begin{equation}
        |a_h^{ij}(x) - a_h^{ij}(y)| \leq C |x-y|^\alpha,
    \end{equation}
    where we used that $Df$ is in $C^{0,\alpha}(B_{r/2})$. The constant $C$ is again independent of $h$. An analogous argument shows that the $\alpha$-Hölder norm of $b_h$ is bounded independently of $h$. Applying the Schauder estimate \cite[Thm.\ 8.32]{GilbargTrudinger}, keeping in mind that the previously discussed bounds do not depend on $h$, we conclude that
    \begin{equation}
        \|v_h\|_{C^{1,\alpha}(B_{r/8})} \leq C
    \end{equation}
    and $C > 0$ is independent of $h$. One can thus again pass to the limit $h \to 0$ to conclude that $f \in C^{2,\alpha}(B_{r/8})$. From here, classical elliptic regularity \cite[Thm.\ 6.17]{GilbargTrudinger} can be used to extrapolate smoothness of $f$: Let us (now classically) rewrite the $p$-harmonicity:
    \begin{equation}
        0 = \partial_i F^i(x,Df(x)) = \partial_{x^i} F^i(x,Df(x)) + \partial_{\xi_j} F^i(x,Df(x)) \partial_{ij}f(x).
    \end{equation}
    Thus, $f \in C^{2,\alpha}(B_{r/8})$ is the solution of an elliptic equation with coefficients of regularity $C^{1,\alpha}(B_{r/8})$, so that $f \in C^{3,\alpha}(B_{r/8})$. Iterating this argument shows that $f \in C^{\infty}(B_{r/8})$. Since the regularity arguments can be translated, we conclude that $f \in C^\infty(\R^n)$, as claimed.
\end{proof}

As a consequence of Theorem \ref{Theorem: smoothness}, we do not need to specify regularity when speaking about $p$-harmonic temporal functions in any of the results below, as they are always automatically smooth.

\section{A Lorentzian Liouville theorem}
\label{section: Liouville}

\subsection{Main theorem and its proof}

Our goal is to prove the following Lorentzian Liouville-type theorem:

\begin{Theorem}[Lorentzian Liouville-type theorem for $p$-harmonic temporal functions]
\label{Theorem: Liouville}
Let $(M,g)$ be a timelike geodesically complete, globally hyperbolic spacetime with compact Cauchy surfaces satisfying the strong energy condition. Let $p < 1$. Suppose $f:M \to \R$ is a temporal function that is $p$-harmonic. Then $f$ is a Cauchy temporal function.
\end{Theorem}

Given the substantial length of the proof, let us give a summary of the key ideas here: The proof proceeds by contradiction. If $f$ remains bounded above along a future inextendible causal curve, we are able to construct a sublevel $\Omega:=\{f<L_*\}$ containing a future-complete timelike ray. Writing $v:=L_*-f$ and $w:=|\nabla f|_g$, divergence-freeness of the vector field $J:=w^{p-2}\nabla f$ gives finite $p$-energy of $f$ on the future part of $\Omega$. The main computational hurdle lies in combining $p$-harmonicity of $f$ with the Raychaudhuri equation for the normalized gradient $\nabla f/w$ and the strong energy condition to obtain an elliptic differential inequality for $G:=w/v^r$, where $r=(n-1)/(n-p)\in(0,1)$. The finite $p$-energy enables us to construct cutoff functions with vanishing weighted energy, which allow us to propagate a positive lower bound for $G$ from an initial compact slab throughout the relevant connected component of the future of $\Omega$, namely the one containing the aforementioned complete timelike ray. Along the latter, the resulting estimate $w\geq\kappa v^r$ forces $v$ to vanish in finite proper time, a contradiction to the fact that the ray is complete and contained in $\Omega$. The analogous argument in the past direction shows that $f$ cannot be bounded along past inextendible causal curves.

\begin{proof}[Proof of Theorem \ref{Theorem: Liouville}]
    Throughout the proof, we set
\begin{equation}
    w:=|\nabla f|_g,\quad N:=\frac{\nabla f}{w},
    \quad q:=1-p, \quad a:=\nabla_NN, \quad A:= D^2H_p(df), \quad \theta:=\mathrm{div}(N).
\end{equation}
Choose a smooth Cauchy temporal function
$\tau:M\to\mathbb R$, which always exists by Bernal--Sánchez \cite{BernalSanchez}, and write
\begin{equation}
    \Sigma_t:=\{\tau = t\}
\end{equation}
for its level sets. The hypersurfaces $\Sigma_t$ are compact, and
$\{t_0 \leq \tau \leq t_1\} = J^+(\Sigma_{t_0}) \cap J^-(\Sigma_{t_1})$ is compact for every finite interval $[a,b]$.
The vector field $T:=\nabla\tau/g(\nabla\tau,\nabla\tau)$ is future timelike and satisfies $T\tau=1$. It is thus easily seen to be complete and any two level sets of $\tau$ are diffeomorphic via the flow of $T$. 

The proof will be established in 13 steps below, we suppose there exists a future inextendible causal curve $\gamma:[0,b) \to M$ such that $\lim_{s \uparrow b} f(\gamma(s)) =: L < + \infty$. We will lead this assumption to a contradiction. Showing that $f$ is unbounded along past inextendible causal curves is analogous.
\medskip\\
\noindent \textbf{STEP 1: Reduction to a sublevel.} 
We consider the function $\mu:\R \to \R$ defined by $\mu(t):=\min_{\Sigma_t} f$. We claim that $\mu$ is strictly increasing. To see this, let $t_1 < t_2$ and $y \in \Sigma_{t_2}$ such that $f(y) = \mu(t_2)$ (such a point exists by compactness). Then there is some $x \in \Sigma_{t_1}$ such that $x \ll y$. Since $f$ is strictly increasing along future causal curves, we get
\begin{equation}
    \mu(t_1) \leq f(x) < f(y) = \mu(t_2).
\end{equation}
Since $\tau$ is Cauchy temporal, $\gamma$ has to meet every $\Sigma_t$ for large $t$ uniquely at a point $z_t \in \Sigma_t \cap \gamma([0,b))$. Thus,
\begin{equation}
    L_*:=\lim_{t \uparrow \infty} \mu(t) \leq \lim_{t \uparrow \infty} f(z_t) = L < + \infty.
\end{equation}
By strict monotonicity, $\mu(t) < L_*$ for all $t \in \R$. Consider the sublevel 
\begin{equation}
\Omega:=\{f < L_*\}
\end{equation}
and set $v:=L_* - f$, thus $v > 0$ on $\Omega$. Now, for every $x \in \Omega \cap \Sigma_t$, we get
\begin{equation}
\label{eq: v leq delta}
    0 < v(x) \leq L_* - \mu(t) \to 0 \quad \text{as} \quad t \to \infty, 
\end{equation}
so we conclude that
\begin{equation}
\label{eq: v to 0 uniformly}
    \lim_{t \to \infty} \sup_{\Omega \cap \Sigma_t} v = 0.
\end{equation}

\noindent \textbf{STEP 2: Replacing $\gamma$ by a timelike ray inside $\Omega$.} Pick a sequence $t_j \to \infty$ as well as points $y_j \in \Sigma_{t_j}$ such that $f(y_j) = \mu(t_j)$. Connect $\Sigma_0$ to $y_j$ by a unit speed maximizing timelike geodesic $\sigma_j:[0,l_j] \to M$. By compactness of $\Sigma_0$, as $j \to \infty$, the $\sigma_j$ converge in $C^\infty_{loc}$ to a timelike unit speed $\Sigma_0$-ray $\sigma:[0,\infty) \to M$. Now for every $s$, we have for large enough $j$
\begin{equation}
    f(\sigma_j(s)) \leq f(y_j) = \mu(t_j) < L_*,
\end{equation}
sending $j \to \infty$ gives
\begin{equation}
    f(\sigma(s)) \leq L_*.
\end{equation}
The above inequality is necessarily strict, since $f(\sigma(s)) = L_*$ would imply $f(\sigma(s')) > L_*$ for any $s' > s$. Thus:
\begin{equation}
    \sigma([0,\infty)) \subseteq \Omega.
\end{equation}
In particular, $v(\sigma(s)) > 0$ for every $s \geq 0$ and $v(\sigma(s)) \to 0$ as $s \to +\infty$ as a consequence of \eqref{eq: v to 0 uniformly}.
\medskip\\
\noindent \textbf{STEP 3: $f$ has finite future $p$-energy in $\Omega$.} We want to show that, for every $t_0 \in \R$, 
\begin{equation}
    \int_{\Omega \cap \{\tau \geq t_0\}} w^p \, dvol_g < + \infty.
\end{equation}
This is the first step in the proof where we use $p$-harmonicity of $f$. We write 
\begin{equation}
J:=w^{p-2} \nabla f.
\end{equation}
Note that $p$-harmonicity of $f$ precisely means that $\mathrm{div}(J) = 0$. Hence,
\begin{equation}
\label{eq: w^p = vJ}
    \mathrm{div}(vJ) = dv(J) = - df(J) = -w^{p-2} df(\nabla f) = -w^p.
\end{equation}
Let us set
\begin{equation}
    F(t):= \int_{\Omega \cap \Sigma_t} v g(J,\nu_t) \, d\mu_t,
\end{equation}
where $\nu_t:=(\nabla \tau/|\nabla \tau|_g)|_{\Sigma_t}$ is the future unit normal of $\Sigma_t$. Since $J$ and $\nu_t$ are future timelike and $v > 0$ on $\Omega$, it follows that $0 \leq F < +\infty$ for all $t \in \R$. Thus, an application of the divergence theorem together with \eqref{eq: w^p = vJ} yields (upon fixing $t_0 < t_1$)
\begin{align}
    -\int_{\Omega \cap \{t_0 \leq \tau \leq t_1\}} w^p  \, dvol_g =\int_{\Omega \cap \{t_0 \leq \tau \leq t_1\}} \mathrm{div}(vJ) \, dvol_g = F(t_1) - F(t_0) \geq -F(t_0).
\end{align}
There is no boundary contribution from $\{f = L_*\}$ since $v = 0$ there. An application of the monotone convergence theorem to send $t_1 \to \infty$ on the left hand side yields
\begin{equation}
    \int_{\Omega \cap \{\tau \geq t_0\}} w^p \, dvol_g \leq F(t_0) < + \infty.
\end{equation}
\noindent \textbf{STEP 4: A useful reexpression of the Raychaudhuri equation for $N$.} Let us write $h:=-g|_{N^\perp}$ for the Riemannian metric induced on the level sets of $f$ and 
\begin{equation}
    K(X,Y):= g(\nabla_X N, Y), \quad X,Y \perp N,
\end{equation}
for the second fundamental form of these level sets. The aim of the current step is to establish the following identity:
\begin{equation}
\label{eq: -div A dw}
    -\mathrm{div}(A(dw)) = w^{p-1} (|K|_h^2 +q|a|_h^2+ \mathrm{Ric}(N,N)).
\end{equation}
First, we note that since $\mathrm{div}(w^{p-1}N) = 0$ by $p$-harmonicity of $f$, we get that
\begin{align}
    0 &= \mathrm{div}(w^{p-1}N) = w^{p-1} \theta + d(w^{p-1})(N) = w^{p-1}\theta +(p-1)w^{p-2} dw (N)\\
    &= w^{p-1} (\theta - qN(\log(w))),
\end{align}
so we conclude the useful identity
\begin{equation}
\label{eq: theta = qNlogw}
    \theta = q N(\log(w)).
\end{equation}
Now using the fact that $\nabla f = wN$ and $g(N,N) = 1$, given a vector field $X \perp N$ we get
\begin{equation}
\label{eq: Hess X N}
    \mathrm{Hess}(f)(X,N) = g(\nabla_X \nabla f, N) = g(X(w)N + w \nabla_X N,N) = X(w).
\end{equation}
In fact, this calculation does not require $X \perp N$, but we are only interested in vector fields tangent to the $f$-levels at the moment. On the other hand,
\begin{align}
\label{eq: Hess N X}
    \mathrm{Hess}(f)(N,X) = g(\nabla_N \nabla f, X) = g(N(w) N + wa,X) = wg(a,X),
\end{align}
where we did use $X \perp N$. By symmetry of the Hessian, \eqref{eq: Hess X N} and \eqref{eq: Hess N X} agree, thus
\begin{equation}
    X(w) = wg(a,X).
\end{equation}
Next, dividing this equation by $w$ implies that $\nabla \log w - a \perp X$. Since $X \perp N$ was arbitrary, we get that $\nabla \log w - a = cN$ for some coefficient function $c$. Since $a = \nabla_NN \perp N$ due to $g(N,N) = 1$, we can calculate $c$ as
\begin{equation}
\label{eq: N-coeff for nabla log w - a}
    c = cg(N,N) = g(\nabla \log w -a,N) = N(\log w).
\end{equation}
On the other hand, taking $X = a$ in $X\log(w) = g(a,X)$ gives
\begin{equation}
\label{eq: a log w = -|a|_h^2}
    a(\log w) = g(a,a) = -|a|_h^2,
\end{equation}
which will be useful later. Let us now recall the expression for $A = D^2 H_p(df)$ from \eqref{eq: D^2 H}, expressed in terms of the conventions in this proof: Given a covector $\alpha$, we get
\begin{equation}
\label{eq: explicit form of A}
    A\alpha = w^{p-2}((2-p)\alpha(N)N - \alpha^\sharp).
\end{equation}
For $\alpha = dw$, 
\begin{align}
    A(dw) &= w^{p-2} ((2-p) N(w)N - \nabla w) = w^{p-1}((2-p)N(\log w)N - \nabla \log w)\\
    &= w^{p-1} ((1-p)N(\log w) N + \underbrace{N(\log w) N - \nabla \log w}_{=-a \text{ by \eqref{eq: N-coeff for nabla log w - a}}})\\
    &= w^{p-1} (\theta N - a), \label{eq: adw = w^p-1 theta N - a}
\end{align}
where the last equality follows from \eqref{eq: theta = qNlogw}, recalling that $q = 1-p$.

Let us now recall the Raychaudhuri equation for not necessarily geodesic unit timelike vector fields (cf.\ Galloway--Ling \cite[Eq.\ (3)]{GallowayLingRemarksCMC} or Hawking--Ellis \cite[Eq.\ 4.26]{HawkingEllis}), which in the case of the vector field $N$ with our conventions reads
\begin{equation}
\label{eq: raychaudhuri for N}
N \theta = \mathrm{div}(a) - |K|_h^2 - \Ric(N,N).
\end{equation}
Now, since $\mathrm{div}(w^{p-1}N) = 0$, we have
\begin{equation}
\label{eq: div w^{p-1} theta N}
    \mathrm{div}(w^{p-1} \theta N) = w^{p-1}N\theta.
\end{equation}
Moreover,
\begin{align}
    \mathrm{div}(w^{p-1}a) &= w^{p-1} \mathrm{div}(a) + (p-1)w^{p-2} dw(a) = w^{p-1} \mathrm{div}(a) + (p-1)\underbrace{w^{p-2} dw(a)}_{=w^{p-1} a(\log w)}\\
    &\overset{\eqref{eq: a log w = -|a|_h^2}}{=} w^{p-1} (\mathrm{div}(a) - q \, a (\log w)). \label{eq: div w^p-1 a = ...}
\end{align}
Thus, combining our previous findings and entering into the Raychaudhuri equation, we find that
\begin{align}
    -\mathrm{div}(A(dw)) &\overset{\eqref{eq: adw = w^p-1 theta N - a}}{=} \mathrm{div}(w^{p-1}a) - \mathrm{div}(w^{p-1}\theta N)\\&\overset{\eqref{eq: div w^p-1 a = ...}}{=} w^{p-1} (\mathrm{div} \, a - q \, a(\log w)) - \mathrm{div}(w^{p-1} \theta N)\\
    &\overset{\eqref{eq: div w^{p-1} theta N}}{=} w^{p-1}(\mathrm{div}\, a - q \, a\log w - N\theta)\\
    &\overset{\eqref{eq: raychaudhuri for N}}{=} w^{p-1}(|K|_h^2 - q \, a \log w + \Ric(N,N)).
\end{align}
\noindent \textbf{STEP 5: $z:=wv^{-1}$ is a subsolution of an elliptic equation.} On $\Omega = \{v > 0\}$, where we recall that $v = L_* - f$, let us set
\begin{equation}
    z:=\frac{w}{v}, \qquad  \mathcal Lu:=-\mathrm{div}(v^2A(du))
\end{equation}
for $u \in C^2(\Omega)$. The aim of this step is to show that $\mathcal L z \geq 0$ on $\Omega$. To this end, let us apply the explicit form of $A$ to $\alpha = dv$ in \eqref{eq: explicit form of A}:
\begin{equation}
    A(dv) = w^{p-2} ((2-p) dv(N)N - \nabla v).
\end{equation}
Since $v = L_* - f$, we have $\nabla v = - \nabla f = -wN$ and $dv(N) = -df(N) = -w$. Thus
\begin{align}
    A(dv) &= w^{p-2}((2-p)(-w)N + wN)\\
    &= (p-1) w^{p-1} N = (p-1) J\\
    &= -qJ.\label{eq: A dv = -qJ}
\end{align}
Thus $\mathrm{div}(A(dv)) = 0$ on $\Omega$ as a consequence of $p$-harmonicity of $f$. Using the product rule, we have $dz = v^{-1}dw - v^{-2}w dv$, hence by linearity of $A$
\begin{equation}
    v^2 A(dz) = v A(dw) - wA(dv).
\end{equation}
Taking $-\mathrm{div}$ on both sides and expanding on the right hand side, we arrive at
\begin{equation}
    \mathcal Lz = -v \,\mathrm{div}(A(dw)) - dv(A(dw)) + w\, \mathrm{div}(A(dv)) + dw(A(dv)).
\end{equation}
The second and fourth terms cancel by symmetry of $A$. The third term vanishes by our observation above that $\mathrm{div}(A(dv)) = -q\, \mathrm{div}(J) = 0$. Thus, by \eqref{eq: -div A dw}, it follows that
\begin{equation}
\label{eq: L z}
    \mathcal L z = -v \, \mathrm{div}(A(dw)) = vw^{p-1} (|K|_h^2 + q|a|_h^2 + \Ric(N,N)) \geq 0,
\end{equation}
where $\Ric(N,N) \geq 0$ uses the strong energy condition and the remaining terms are manifestly nonnegative.
\medskip\\
\noindent \textbf{STEP 6: Computing $\mathcal L \log z$.} Given that $d \log z = z^{-1} dz$, we have
\begin{align}
    \mathcal L \log z &= -\mathrm{div}(v^2 A(d\log z)) = -\mathrm{div}\left(\frac{v^2}{z} A(dz)\right)\\
    &= \frac{\mathcal{L}z}{z} + \frac{v^2}{z^2} dz(A(dz)) \\
    &= \frac{\mathcal{L}z}{z} + v^2 d\log z(A(d\log z)). \label{eq: L log z = I + II}
\end{align}
Let us compute the two terms on the right hand side. For the first term, using \eqref{eq: L z} and the definition of $z = wv^{-1}$, we have
\begin{equation}
\label{eq: Lz /z}
    \frac{\mathcal L z}{z} = v^2 w^{p-2}(|K|_h^2 + q|a|_h^2 + \Ric(N,N)).
\end{equation}
For the second term, let us set $x:=N \log w$ for simplicity. We know from \eqref{eq: N-coeff for nabla log w - a} that
\begin{equation}
    \nabla \log w = xN + a.
\end{equation}
Moreover, since $\nabla v = -wN$, we have
\begin{equation}
    \nabla \log v = v^{-1} \nabla v = -\frac{w}{v} N = -zN.
\end{equation}
Combining these two, we arrive at
\begin{equation}
    \nabla \log z = \nabla (\log w - \log v) = (x + z)N + a.
\end{equation}
Consequently, $d \log z(N) = x + z$. Using again the explicit form of $A$ from \eqref{eq: explicit form of A} applied to $d\log z$, we get
\begin{align}
    A(d\log z) &= w^{p-2}((2-p) d\log z(N) N - \nabla \log z)\\
    &= w^{p-2}((2-p) (x+z)N - (x+z)N - a)\\
    &= w^{p-2}(q(x+z)N - a).
\end{align}
Combining these expressions for $\nabla \log z$ and $A(d\log z)$, we get
\begin{align}
    d\log z (A(d\log z)) &= g(A(d\log z), \nabla \log z) = g(w^{p-2}(q(x+z)N - a), (x+z)N + a)\\
    &= w^{p-2}(q(x+z)^2 - g(a,a))\\
    &= w^{p-2}(q(x+z)^2 + |a|_h^2), \label{eq: dlog z A d log z}
\end{align}
where we used that $a \perp N$ (which follows from $g(N,N) = 1$) and $g(a,a) = -|a|_h^2$ (cf.\ \eqref{eq: a log w = -|a|_h^2}). Reinserting the expressions from \eqref{eq: Lz /z} and \eqref{eq: dlog z A d log z} back into \eqref{eq: L log z = I + II}, we conclude that
\begin{equation}
\label{eq: L log z}
    \mathcal L \log z =v^2 w^{p-2} (|K|_h^2 + (1+q)|a|_h^2 + \Ric(N,N) + q(x+z)^2).
\end{equation}
\noindent \textbf{STEP 7: Computing $\mathcal L \log v$.} Since $v^2 A(d \log v) = v A(dv) = -qvJ$ due to \eqref{eq: A dv = -qJ}, we get
\begin{equation}
    \mathcal{L }\log v = q\, \mathrm{div}(vJ) = q \, dv(J) = -qw^p,
\end{equation}
where we used $\mathrm{div}(J) = 0$ and \eqref{eq: w^p = vJ}. Thus,
\begin{equation}
    \mathcal L \log v = -qw^p = -qv^2 w^{p-2} z^2,
\end{equation}
where the definition $z = w/v$ was used.
\medskip\\
\noindent \textbf{STEP 8: Completing the square.}
We now interpolate between $\log z$ and $\log v$, with the intention of finding $r\in(0,1)$ such that $\psi=\log(w/v^r)$ satisfies $\mathcal L\psi\geq v^2d\psi(A(d\psi))$, ensuring that $G = e^\psi$ satisfies $\mathcal L G \geq 0$. The exponent $r \in (0,1)$ will turn a (yet to be established) positive lower bound for $G$ into a finite-time contradiction along the complete timelike ray. We thus define
\begin{equation}
    \psi := \log z + (1-r) \log v = \log w - \log v + (1-r)\log v = \log\left( \frac{w}{v^r}\right),
\end{equation}
where $r$ is to be determined. Let us compute $\mathcal L \psi$, using linearity as well as the expressions for $\mathcal L \log z$ and $\mathcal L \log v$ obtained in Steps 6 and 7:
\begin{align}
    \mathcal L \psi &= \mathcal L \log z + (1-r)\mathcal L \log v \\
    &= v^2 w^{p-2}(|K|_h^2 + (1+q)|a|_h^2 + \Ric(N,N) + q((x+z)^2 - (1-r)z^2))\\
    &=v^2 w^{p-2}(|K|_h^2 + (1+q)|a|_h^2 + \Ric(N,N) + q(x^2 + 2xz + rz^2)).\label{eq: L psi}
\end{align}
We have $\mathrm{Tr}_h K = -\theta = -qx$, where the first equality is due to our sign conventions and the second follows from \eqref{eq: theta = qNlogw} together with the definition $x = N\log w$. Thus, by Cauchy-Schwarz,
\begin{equation}
    q^2 x^2 = \theta^2 \leq (n-1) |K|_h^2,
\end{equation}
so that
\begin{equation}
    |K|_h^2 \geq \frac{q^2}{n-1} x^2.
\end{equation}
Using this estimate in \eqref{eq: L psi}, we get
\begin{equation}
    \mathcal L \psi \geq v^2 w^{p-2} \bigg((1+q)|a|_h^2 + \Ric(N,N) + q\bigg((1+\frac{q}{n-1})x^2 + 2xz + rz^2\bigg)\bigg).
\end{equation}
We now make the choice $r:=(n-1)/(n-p) \in (0,1)$, so that $1/r = 1 + q/(n-1)$, in order to complete the square:
\begin{equation}
    q\bigg(\bigg(1 + \frac{q}{n-1}\bigg)x^2 + 2xz + rz^2 \bigg) = \frac{q}{r} \bigg(x + rz\bigg)^2.
\end{equation}
Finally, since $\psi = \log w - r\log v$, we get
\begin{equation}
\label{eq: Npsi}
    N \psi = x - r N\log v = x + rv^{-1}w = x+ rz.
\end{equation}
All in all,
\begin{equation}
\label{eq: L psi inequality}
    \mathcal L \psi \geq v^2 w^{p-2} \bigg((1+q)|a|_h^2 + \Ric(N,N) + \frac{q}{r} (N\psi)^2\bigg)
\end{equation}
\noindent \textbf{STEP 9: Exponentiating $\psi$.}
We set
\begin{equation}
    G:=e^\psi = \frac{w}{v^r}, \quad \text{ where } \quad r = \frac{n-1}{n-p} \in (0,1).
\end{equation}
Our goal of Step 9 is to show that $\mathcal L G \geq 0$. We start with the following computation, using $dG = G d\psi$:
\begin{align}
    \mathcal L G &= -\mathrm{div}(v^2 A(dG)) = -\mathrm{div}(v^2 G A(d\psi))\\
    &= G(\mathcal L\psi - v^2 d\psi A(d\psi)). \label{eq: expression for L G}
\end{align}
We want to calculate the second term. To this end, observe (using $\nabla \log v = - zN$) that
\begin{equation}
\label{eq: nabla psi}
    \nabla \psi = \nabla \log w - r\nabla \log v  \overset{\eqref{eq: N-coeff for nabla log w - a}}{=} (x+rz)N + a \overset{\eqref{eq: Npsi}}{=} (N \psi)N + a.
\end{equation}
Using again the explicit formula for $A$ from \eqref{eq: explicit form of A}, applied to $d\psi$, we get
\begin{equation}
\label{eq: A dpsi}
    A(d\psi) = w^{p-2}((2-p) N\psi N - \nabla \psi) \overset{\eqref{eq: nabla psi}}{=} w^{p-2}(q(N\psi)N - a).
\end{equation}
Thus, using \eqref{eq: A dpsi} and \eqref{eq: nabla psi}, as well as $g(a,a) = -|a|_h^2$, $g(N,N) = 1$ and $g(N,a) = 0$, we get
\begin{align}
\label{eq: d psi A d psi}
    d\psi A(d\psi) &= g(A(d\psi),\nabla \psi) = g(w^{p-2}(q(N\psi)N - a), N(\psi)N + a) \\
    &=w^{p-2} (q (N\psi)^2 + |a|_h^2).
\end{align}
Using now \eqref{eq: L psi inequality} and \eqref{eq: d psi A d psi} in the expression \eqref{eq: expression for L G}, we get
\begin{align}
    \mathcal L G &\geq G v^2 w^{p-2}\bigg((1+q)|a|_h^2 + \Ric(N,N) + \frac{q}{r} (N\psi)^2 - q(N\psi)^2 - |a|_h^2\bigg)\\
    &= Gv^2 w^{p-2}\bigg(q |a|_h^2 + \Ric(N,N) + \frac{q^2}{n-1} (N\psi)^2\bigg) \geq 0.
\end{align}
where we used $q/r - q = q^2/(n-1)$. The nonnegativity uses the strong energy condition as well as the fact that $q > 0$.
\medskip\\
\noindent \textbf{STEP 10: Constructing cutoffs whose $\mathcal L$-energy tends to $0$.}
Let $E$ be the connected component of $\Omega \cap \{\tau > t_0\}$ (recall the use of the latter from Step 3) containing the complete timelike $\Sigma_0$-ray $\sigma$ up to a compact segment. Recall that, as a consequence of Step 3,
\begin{equation}
    \int_E w^p \, dvol_g < +\infty.
\end{equation}
Note that, since $\mu(t)=\min_{\Sigma_t}f$ is strictly increasing (this was established in Step 1), we get that
\begin{equation}
    0 < v \leq L_* - \mu(t_0) =:\delta_0 \quad \text{on }E.
\end{equation}
Thus,
\begin{equation}
    \rho:=\log \frac{\delta_0}{v} \geq 0 \quad \text{ on }E.
\end{equation}
Since $d\rho = -v^{-1}dv$, $A(dv) = -qJ$ by \eqref{eq: A dv = -qJ} and $dv(J) = -w^p$ by \eqref{eq: w^p = vJ}, we get
\begin{equation}
\label{eq: v2drhoAdrho}
    v^2 d\rho (A(d\rho)) = dv (A(dv)) = -q \, dv(J) = qw^p.
\end{equation}
For $R > 0$, we define the following locally Lipschitz function:
\begin{equation}
    \eta_R:=\bigg(1 - \frac{\rho}{R}\bigg)_+.
\end{equation}
By definition, $0 \leq \eta_R \leq 1$, and $\eta_R \uparrow 1$ pointwise as $R \to \infty$. Since $\eta_R = 0$ whenever $R \leq \rho$, it follows using \eqref{eq: v2drhoAdrho} that
\begin{align}
    \int_E v^2 d \eta_R (A(d\eta_R)) \, dvol_g &= \frac{1}{R^2}\int_{E \cap \{\rho < R\}} v^2d\rho (A(d\rho)) \, dvol_g\\
    &\leq \frac{q}{R^2} \int_E w^p \, dvol_g.
\end{align}
By Step 3, the right hand side converges to $0$ as $R \to \infty$, thus
\begin{equation}
\label{eq: cutoff energy estimate}
    \int_E v^2 d \eta_R (A(d\eta_R)) \, dvol_g \to 0 \quad \text{as} \quad R \to \infty.
\end{equation}
\noindent \textbf{STEP 11: Replacing $\eta_R$ by $\varphi_R$ which is compactly supported in $E$.}
We will now find $\kappa > 0$ such that
\begin{equation}
    \varphi_R:=(\kappa - G)_+ \eta_R^2 
\end{equation}
has compact support in the open set $E$. To this end, we fix $\varepsilon > 0$. Since $\tau^{-1}([t_0,t_0 + \varepsilon]) = J^+(\Sigma_{t_0}) \cap J^-(\Sigma_{t_0 + \varepsilon}) \subseteq M$ is compact by global hyperbolicity, $w$ assumes a minimum $w_\varepsilon > 0$ on it. Since $0 < v \leq \delta_0$ on $E$, we have
\begin{equation}
    G = \frac{w}{v^r} \geq \frac{w_\varepsilon}{\delta_0^r} \quad \text{on} \quad E \cap \tau^{-1}((t_0,t_0+\varepsilon)).
\end{equation}
We choose $\kappa:=w_\varepsilon/2\delta^r_0 >0$ and set $\zeta:=(\kappa - G)_+$. Then $0 \leq \zeta \leq \kappa$, and $\zeta = 0$ on the set $E \cap \tau^{-1}((t_0,t_0 + \varepsilon))$. Next, let us recall that $\eta_R > 0$ if and only if $\rho < R$, which by definition of $\rho = \log \delta_0/v$ can be rewritten as $v > \delta_0 e^{-R}$. Since $v \leq \delta \circ \tau$ from \eqref{eq: v leq delta}, where $\delta(t):=L_*-\mu(t) \to 0$ as $t \to \infty$, for fixed $R > 0$ we can find $T_R > 0$ such that
\begin{equation}
    v < \delta_0 e^{-R} \quad \text{on} \quad E \cap \{\tau \geq T_R\}.
\end{equation}
In particular, $\eta_R = 0$ and thus $\varphi_R = 0$ on $E \cap \{\tau \geq T_R\}$. Hence $\mathrm{supp}(\varphi_R) \subseteq E \cap \{t_0 + \varepsilon \leq \tau \leq T_R, \, v \geq \delta_0 e^{-R}\}$ and we conclude that, for each $R > 0$, the support of the Lipschitz function $\varphi_R$ is compact and contained in $E$.
\medskip\\
\noindent \textbf{STEP 12: Using $\mathcal L G \geq 0$ to show $G \geq \kappa$.} We will now test the already established inequality $\mathcal L G \geq 0$ (cf.\ Step 9) against the functions $\varphi_R$ from Step 11 to conclude that $G \geq \kappa$. To this end, we write
\begin{equation}
    \mathcal B(\alpha,\beta):=v^2 \alpha(A(\beta))
\end{equation}
for covectors $\alpha,\beta$. By symmetry and positive definiteness of $A$, $\mathcal B$ is a symmetric positive definite bilinear form on covectors. Using now $\mathcal L G \geq 0$ and integrating by parts, we get
\begin{equation}
\label{eq: 0 leq int B}
    0 \leq \int_E \varphi_R \mathcal L G \, dvol_g = \int_E \mathcal B(d\varphi_R,dG)\, dvol_g.
\end{equation}
By the product rule, we have
\begin{equation}
    d\varphi_R = \eta_R^2 d\zeta + 2\eta_R \zeta d\eta_R
\end{equation}
almost everywhere. Note that $d\zeta = -dG$ a.e.\ on $\{G < \kappa\}$ and $\zeta = 0$ elsewhere, in particular $d\zeta = 0$. Thus, continuing from \eqref{eq: 0 leq int B}, we have
\begin{align}
    0 &\leq \int_E \mathcal{B}(\eta_R^2 d\zeta + 2\eta_R\zeta d\eta_R,dG) \, dvol_g\\
    &= -\underbrace{\int_E \eta_R^2\mathcal B(d\zeta,d\zeta) \, dvol_g}_{=:I_R} - 2 \int_E \eta_R \zeta\mathcal B(d\eta_R,d\zeta) \, dvol_g.
\end{align}
This immediately implies that
\begin{equation}
\label{eq: I_R estimate}
    0 \leq I_R \leq 2 \int_E\eta_R \zeta |\mathcal B(d\eta_R,d\zeta)| \, dvol_g.
\end{equation}
Using Cauchy-Schwarz and Young's inequality for the positive definite form $\mathcal B$, we get that a.e.
\begin{align}
    2\eta_R \zeta |\mathcal B(d\eta_R,d\zeta)| &\leq 2\eta_R \zeta \sqrt{\mathcal B(d\eta_R,d\eta_R)} \sqrt{\mathcal B(d\zeta,d\zeta)} \\
    &\leq \frac{1}{2}\eta_R^2 \mathcal B(d\zeta,d\zeta) + 2 \zeta^2 \mathcal B(d\eta_R,d\eta_R).
\end{align}
Integrating, using $0 \leq \zeta \leq \kappa$ and the estimate \eqref{eq: I_R estimate}, we get
\begin{equation}
    I_R \leq \frac{1}{2}I_R + 2\kappa^2 \int_E \mathcal B(d\eta_R,d\eta_R) \, dvol_g.
\end{equation}
Thus,
\begin{equation}
    I_R \leq 4\kappa^2 \int_E \mathcal B(d\eta_R,d\eta_R) \, dvol_g \to 0 \quad \text{as} \quad R \to \infty,
\end{equation}
where we used \eqref{eq: cutoff energy estimate} from Step 10 for the limit $R \to \infty$. Since $\eta_R \uparrow 1$ pointwise as $R \to \infty$, monotone convergence implies
\begin{equation}
    0 = \lim_{R \to \infty} I_R = \int_E \mathcal B(d\zeta,d\zeta) \, dvol_g.
\end{equation}
By positive definiteness of $\mathcal B$, we conclude that $d\zeta = 0$ a.e.\ on $E$. Since $\zeta = 0$ on $E \cap \tau^{-1}((t_0,t_0 + \varepsilon))$ (cf.\ its construction from Step 11), we arrive at $\zeta \equiv 0$ on $E$. Thus,
\begin{equation}
    G \geq \kappa \quad \text{on} \quad E,
\end{equation}
or, equivalently,
\begin{equation}
    w \geq \kappa v^r \quad \text{on} \quad E.
\end{equation}
\noindent \textbf{STEP 13: Leading $w \geq \kappa v^r$ to a contradiction.}
Recall that $\sigma:[0,\infty) \to M$ is a future complete unit speed timelike $\Sigma_0$-ray. By definition of $E$, there is $s_1 \geq 0$ such that $\sigma|_{[s_1,\infty)} \subseteq E$. We set $y(s):=v(\sigma(s)) > 0$, $s \geq s_1$. Using $g(N,\dot \sigma) \geq 1$, which follows from the reverse Cauchy-Schwarz inequality, we get
\begin{align}
    y'(s) &= dv(\dot \sigma(s)) = -w(\sigma(s))g(N(\sigma(s)),\dot \sigma(s)) \\
    &\leq -w(\sigma(s)) \leq -\kappa v(\sigma(s))^r = -\kappa y(s)^r.
\end{align}
Multiplying by the positive quantity $(1-r)y(s)^{-r}$ yields
\begin{equation}
    \frac{d}{ds} (y^{1-r}) \leq -(1-r)\kappa.
\end{equation}
Integrating from $s_1$ to $s$ gives
\begin{equation}
    0 < y(s)^{1-r} \leq y(s_1)^{1-r} - (1-r)\kappa(s - s_1).
\end{equation}
The right hand side becomes zero at
\begin{equation}
    s = s_1 + \frac{y(s_1)^{1-r}}{(1-r)\kappa},
\end{equation}
which contradicts the fact that $y(s) > 0$ for all $s \geq s_1$.
\end{proof}

\subsection{Necessity of the assumptions}
\label{subsection: necessity of assumptions}

In Theorem \ref{Theorem: Liouville}, all of the assumptions, which are precisely the ones in Bartnik's conjecture, are necessary, as we will now show. This underscores the fact that Liouville-type unboundedness of $p$-harmonic temporal functions is intimately connected to Bartnik's conjecture.

\begin{example}[Timelike geodesic completeness is necessary]
Let $(M,g) = ((a,b) \times \Sigma, dt^2 - h)$, where $(\Sigma,h)$ is a Riemannian manifold and either $a > -\infty$ or $b < + \infty$. Then $f((t,x)):=t$ is temporal and $p$-harmonic for every $p$, but it is clearly not Cauchy temporal.
\end{example}

\begin{example}[Compactness of Cauchy surfaces is necessary]
Consider $2$-dimensional Minkowski spacetime $(\mathbb{R}^2,g = dt^2 - dx^2)$. Fix $p<1$ and set $q:=1-p>0$. Define
\begin{equation}
    X(s):=\frac{(1+q)\sinh(2s)+2(1-q)s}{4q}.
\end{equation}
Since
\begin{equation}
    X'(s)=\frac{1+(1+q)\sinh^2s}{q}>0
\end{equation}
and $X(s)\to\pm\infty$ as $s\to\pm\infty$, the map $X:\R \to \R$ is a
smooth diffeomorphism. Thus
\begin{equation}
    f(t,X(s))
    :=-e^{-t-\frac{1+q}{2q}\sinh^2s}\cosh s
\end{equation}
defines a smooth function on all of Minkowski spacetime. Writing $s=X^{-1}(x)$ and
$w=e^{-t-\frac{1+q}{2q}\sinh^2s} > 0$, direct differentiation gives
\begin{equation}
    \partial_t f=w\cosh s,
    \qquad
    \partial_x f=w\sinh s.
\end{equation}
Hence $\nabla f$ is future timelike and $|\nabla f|_g=w>0$.
Moreover,
\begin{equation}
    J:=w^{p-2}\nabla f
      =w^{-q}\bigl(\cosh s\,\partial_t-\sinh s\,\partial_x\bigr).
\end{equation}
Since
\begin{equation}
    \partial_t J^t=q w^{-q}\cosh s,
    \qquad
    \partial_x J^x=-q w^{-q}\cosh s,
\end{equation}
we obtain $\square_p f=-\operatorname{div}J=0$. However, $f < 0$ everywhere, thus it is not Cauchy temporal.
\end{example}

\begin{example}[The strong energy condition is necessary]
Consider $M:=\mathbb{R}\times S^1$ with the warped product metric
\begin{equation}
    g:=dt^2-\frac{1}{1+t^2}\,dx^2,
    \qquad
    f(t,x):=\arctan t.
\end{equation}
The slices $\{t\}\times S^1$ are compact Cauchy hypersurfaces.
Writing $a(t)=(1+t^2)^{-1/2}$, we have
\begin{equation}
    \nabla f=\frac{1}{1+t^2}\,\partial_t,
    \qquad
    w:=|\nabla f|_g=\frac{1}{1+t^2}=a^2.
\end{equation}
Thus $f$ is smooth and temporal. For $p=\frac12$, we obtain
\begin{equation}
    \square_pf
    =-\frac{1}{a}\frac{d}{dt}\bigl(a\,w^{-1/2}\bigr)
    =0.
\end{equation}
However, $f$ is bounded and therefore not Cauchy temporal. The spacetime $(M,g)$ is easily seen to be timelike geodesically complete, however
\begin{equation}
    \operatorname{Ric}(\partial_t,\partial_t)
    =-\frac{a''}{a}
    =\frac{1-2t^2}{(1+t^2)^2}<0
    \qquad\text{for }|t|>\frac{1}{\sqrt{2}},
\end{equation}
so the strong energy condition fails.
\end{example}

\section{Cosmological splitting via $p$-harmonic Cauchy temporal functions}
\label{section: cosmological rigidity}

The purpose of this section is to show that cosmological spacetimes which contain a $p$-harmonic Cauchy temporal function necessarily split as a product. Let us emphasize that this result is not a consequence of the known splitting under the existence of a CMC Cauchy surface \cite{Bartnik1988remarks}, since the norm of the gradient of a $p$-harmonic function is not a priori constant. To establish the latter, we will make use the $p$-Bochner--Ohta identity from Theorem \ref{Theorem: p Bochner Ohta}.

Before coming to this, we include a lemma on the derivative with respect to the level set parameter of integrals on level sets of a Cauchy temporal function, which follows from the more general results of Gerhardt \cite{Gerhardt} on the evolution of metrics.

\begin{Lemma}[First variation of level set integrals]
\label{Lemma: first variation of level set}
Let $(M,g)$ be a globally hyperbolic spacetime with compact Cauchy surfaces. Let $f:M \to \R$ be a smooth Cauchy temporal function. Denote by $\Sigma_s:=\{f =s \}$ its $s$-level set and by $\mu_s$ the corresponding hypersurface measure. Then for every $\varphi \in C^\infty(M)$ and every $s \in \R$, 
\begin{equation}
    \frac{d}{ds} \int_{\Sigma_s} \varphi \, d\mu_s = \int_{\Sigma_s} \bigg[\frac{df(\nabla \varphi)}{|\nabla f|^2_g} + \frac{\varphi}{|\nabla f|_g} \mathrm{div}\bigg(\frac{\nabla f}{|\nabla f|_g}\bigg) \, \bigg]d\mu_s.
\end{equation}
\end{Lemma}
\begin{proof}
Let us set $w:=|\nabla f|_g$, $N:=\nabla f/w$ and $V:=N/w$. Then $df(V) = 1$ and $V$ is consequently a complete vector field: Indeed, if $\gamma:(a,b) \to M$ is a flow line of $V$, then it is an inextendible causal curve and satisfies
\begin{equation}
    (f \circ \gamma)'(s) = df(V(\gamma(s))) = 1,
\end{equation}
so that $f(\gamma(s)) = s + f(\gamma(0))$. Since $f$ is a Cauchy temporal function, it follows that the left hand side goes to $+\infty$ for $s \uparrow b$, thus $b = + \infty$. Analogously, $a = -\infty$ and $\gamma$ is defined on all of $\R$. From here, it is easily seen that the flow $F_s:\Sigma_0 \to \Sigma_s$ mapping $x\mapsto \gamma_x(s)$, where $\gamma_x$ is the unique flow line of $V$ with $\gamma_x(0) = x$, is a diffeomorphism. Thus, we may write
\begin{equation}
    \int_{\Sigma_s} \varphi \, d\mu_s = \int_{\Sigma_0} (\varphi \circ F_s) \, d(F_s^*\mu_s).
\end{equation}
Writing $\nu_s := F_s^* \mu_s = \sqrt{\det q_s} \,|dx^1 \wedge \dots dx^{n-1}|$ in local coordinates, where $q_s$ denotes the $F_s$-pullback of the Riemannian metric $-g|_{\Sigma_s}$ on $\Sigma_s$ onto $\Sigma_0$, the claim now follows from the product rule, Jacobi's formula to differentiate the determinant and using \cite[Lem.\ 3.2]{Gerhardt} to calculate $\partial_s q_s$.
\end{proof}

\begin{Theorem}[Cosmological splitting via $p$-harmonic Cauchy temporal functions]
\label{theorem: cosmological splitting p harmonicity}
Let $(M,g)$ be a globally hyperbolic spacetime with compact Cauchy surfaces satisfying the strong energy condition. Suppose there exists a smooth Cauchy temporal function $f:M \to \R$ which is $p$-harmonic for some $p < 1$, i.e., it satisfies
\begin{equation}
    \square_p f=-\mathrm{div}(|\nabla f|_g^{p-2} \nabla f) = 0.
\end{equation}
Then $(M,g)$ splits isometrically as $(\R \times \Sigma, dt^2 - \tilde g)$.
\end{Theorem}
\begin{proof}
Fix $p < 1$ for which $f$ is $p$-harmonic. Let us recall the $p$-Bochner--Ohta identity from Theorem \ref{Theorem: p Bochner Ohta}, which states that
\begin{equation}
    -\mathrm{div}(D^2H_p|_{df} [d(u_p(|df|_{g^*}))]) - DH_p|_{df} [d \square_p f] = |\mathrm{Hess} f|_{D^2H_p|_{df}}^2 + \mathrm{Ric}(DH_p|_{df},DH_p|_{df}).
\end{equation}
Here, $u_p(t)=p^{-1}t^p$ if $0 \neq p < 1$ and $u_0(t)=\log(t)$. The second term on the left hand side vanishes due to $p$-harmonicity of $f$. We will use the following notation throughout the proof: We write $A_p:=D^2H_p|_{df}$, $\Sigma_s = \{f = s\}$, $w:=|\nabla f|_g$, $N:=|\nabla f|_g^{-1} \nabla f$, as well as $\theta :=\mathrm{div}N$. It follows that
\begin{equation}
    - \mathrm{div}\left(A_p [d (u_p(w))]\right) = Q:=|\mathrm{Hess} \, f|^2_{D^2H_p|_{df}}  + w^{2p-2}\mathrm{Ric}(N,N) \geq 0,
\end{equation}
where we have expanded the explicit form of $H_p$ and its derivatives, and the last inequality is a consequence of the strong energy condition. Note that $Q = 0$ implies $\mathrm{Hess} \, f = 0$. Since
\begin{align}
    0 &= -\square_p f = \mathrm{div}(|\nabla f|_g^{p-2} \nabla f) = \mathrm{div}(|\nabla f|_g^{p-1} N) = d (w^{p-1})(N) + w^{p-1} \mathrm{div}N\\
    &= (p-1) w^{p-2} N(w) + w^{p-1} \theta,
\end{align}
we conclude the useful formula
\begin{equation}
\label{eq: theta}
    \theta = (1-p)\frac{N(w)}{w}.
\end{equation}
Further,
\begin{align}
    Z_p:=A_p d(u_p(w)) &= A_p w^{p-1} dw = w^{p-2}[(2-p) N \otimes N - g^*] [w^{p-1} dw] \\
    &= w^{2p-3} [(2-p) N(w) N - \nabla w]\\
    &= w^{2p-3}[(1-p) N(w)N + \nabla^{\Sigma_s} w],
\end{align}
where we used $\nabla w= N(w) N - \nabla^{\Sigma_s} w$ on $\Sigma_s:=\{f =s\}$. The minus sign is because of the signature choice of $g$, by $\nabla^{\Sigma_s} w$ we of course mean the gradient with respect to the \emph{negative of the restriction} of $g$ to $\Sigma_s$, which is a Riemannian metric. Thus, on $\Sigma_s$,
\begin{equation}
    g(Z_p,N) = w^{2p-3} (1-p) N(w) = w^{2p-2} \theta.
\end{equation}
We define $B_p:\R \to \R$ via
\begin{equation}
\label{def: B_p}
    B_p(s):=\int_{\Sigma_s} g(Z_p,N) \, d\mu_s = \int_{\Sigma_s} w^{2p-2} \theta \, d\mu_s.
\end{equation}
An application of the divergence theorem yields (with $\Omega_{a,b}:=\{a \leq f \leq b\}$)
\begin{align}
\label{eq: int Q = diff of Bs}
    \int_{\Omega_{a,b}} Q \, dvol_g &= -\int_{\Omega_{a,b}} \mathrm{div} \, Z_p \, dvol_g = \int_{\Sigma_a} g(Z_p,N) \, d\mu_a - \int_{\Sigma_b} g(Z_p,N) \, d\mu_b \nonumber\\
    &= B_p(a) - B_p(b).
\end{align}
Since $Q \geq 0$, it follows that $B_p$ is nonincreasing. By the coarea formula, the integral of $Q$ over the region $\Omega_{a,b}$ can be expressed as
\begin{equation}
    \int_{\Omega_{a,b}} Q \, \mathrm{dvol}_g = \int_a^b \int_{\Sigma_s} \frac{Q}{w} \, d\mu_s \, ds.
\end{equation}
Thus, it follows from \eqref{eq: int Q = diff of Bs} that
\begin{equation}
\label{eq: derivative of B_p}
    B_p'(s) = -\int_{\Sigma_s} \frac{Q}{w} \, d\mu_s \leq 0.
\end{equation}
We define now the function $C_p:\R \to \R$ by
\begin{equation}
\label{eq: def C_p}
    C_p(s):= \int_{\Sigma_s} w^{2p-1} \, \mathrm d \mu_s.
\end{equation}
By compactness of $\Sigma_s$, we have $0 < C_p(s)<+\infty$ for all $s \in \R$. By Lemma \ref{Lemma: first variation of level set}, we compute
\begin{align}
    C'_p(s) &= \int_{\Sigma_s} \frac{N(w^{2p-1})}{w} + \frac{\theta}{w} w^{2p-1} \, \mathrm d\mu_s\\
    &= \int_{\Sigma_s} (2p-1) w^{2p-3} N(w) + w^{2p-2} \theta \, \mathrm d\mu_s \\
    &\overset{\eqref{eq: theta}}{=} \int_{\Sigma_s} (2p -1 + 1 - p) w^{2p-3} N(w) \, \mathrm d\mu_s\\
    &\overset{\eqref{eq: theta}}{=} \frac{p}{1-p} \int_{\Sigma_s} w^{2p-2} \theta \, \mathrm d\mu_s\\
    &= \frac{p}{1-p} B_p(s). \label{eq: C_p' = p/q B_p}
\end{align}
Put $q:=1-p>0$. At any point, choose a Lorentzian orthonormal frame
$e_0=N,e_1,\ldots,e_{n-1}$ and write
$H_{ab}:=\operatorname{Hess}f(e_a,e_b)$.
Since $A_p=w^{p-2}\operatorname{diag}(q,1,\ldots,1)$ and
$H_{ab}=H_{ba}$, we have
\begin{equation}
    |\operatorname{Hess}f|_{A_p}^2
    =w^{2p-4}\left(
        q^2H_{00}^2
        +2q\sum_{i=1}^{n-1}H_{0i}^2
        +\sum_{i,j=1}^{n-1}H_{ij}^2
    \right).
\end{equation}
Moreover, $\nabla f=wN$ and $g(N,N)=1$ imply
\begin{equation}
    H_{00}
    =g(\nabla_N(wN),N)
    =N(w)+\frac{w}{2}N(g(N,N))
    =N(w).
\end{equation}
Thus, since $q>0$ and $\operatorname{Ric}(N,N)\ge0$,
\begin{equation}
\label{eq: bound on Q}
    Q=|\operatorname{Hess}f|_{A_p}^2
      +w^{2p-2}\operatorname{Ric}(N,N)
    \ge q^2w^{2p-4}(N(w))^2.
\end{equation}
Using
\begin{equation}
    B_p(s)=q\int_{\Sigma_s}w^{2p-3}N(w)\,d\mu_s,
    \qquad
    C_p(s)=\int_{\Sigma_s}w^{2p-1}\,d\mu_s,
\end{equation}
where the first follows from the definition of $B_p$ \eqref{def: B_p} together with the identity \eqref{eq: theta}, we apply the Cauchy--Schwarz inequality to the factors
$w^{p-\frac12}$ and $w^{p-\frac52}N(w)$ to obtain
    \begin{align}
        B_p(s)^2
        &=q^2\left(
            \int_{\Sigma_s}
            w^{p-\frac12}w^{p-\frac52}N(w)\,d\mu_s
        \right)^2\\
        &\le q^2
        \left(\int_{\Sigma_s}w^{2p-1}\,d\mu_s\right)
        \left(\int_{\Sigma_s}w^{2p-5}(N(w))^2\,d\mu_s\right)\\
        &=q^2C_p(s)
        \int_{\Sigma_s}w^{2p-5}(N(w))^2\,d\mu_s.
    \end{align}
Dividing \eqref{eq: bound on Q} by $w>0$ and integrating,
\eqref{eq: derivative of B_p} gives
\begin{equation}
    q^2\int_{\Sigma_s}w^{2p-5}(N(w))^2\,d\mu_s
    \le\int_{\Sigma_s}\frac{Q}{w}\,d\mu_s
    =-B_p'(s).
\end{equation}
Since $C_p(s)>0$, combining these inequalities yields
\begin{equation}
\label{eq: B_p squared ineq}
    B_p(s)^2\le -C_p(s)B_p'(s).
\end{equation}
Define the positive function $D_s:\R \to (0,\infty)$ (recalling that $C_p > 0$)
\begin{equation}
    D(s):=\exp\left(
        \int_0^s\frac{B_p(\sigma)}{qC_p(\sigma)}\,d\sigma
    \right).
\end{equation}
Using $C_p'=(p/q)B_p$ by \eqref{eq: C_p' = p/q B_p}, we compute using \eqref{eq: B_p squared ineq} 
\begin{equation}
    D''(s)
    =\frac{D(s)}{qC_p(s)^2}
    \bigl(C_p(s)B_p'(s)+B_p(s)^2\bigr)
    \le 0.
\end{equation}
Thus $D$ is positive and concave on $\mathbb{R}$, hence constant.
Its definition then implies $B_p=0$ on $\mathbb{R}$ since $qC_p > 0$. Since
\begin{equation}
    \int_{\Omega_{a,b}} Q \, \mathrm{dvol}_g = B_p(a) - B_p(b) = 0
\end{equation}
for every $a < b$ and $Q \geq 0$, we conclude that indeed $Q = 0$ on all of $M$, from which $\mathrm{Hess} \, f = 0$ follows. In particular, $w = |\nabla f|_g \equiv c > 0$ is constant on $M$.

Upon replacing $f$ with $f/c$, we may assume that $c = 1$. We now show that $\nabla f$ is a complete vector field. Indeed, let $\gamma:(a,b) \to M$ be a maximal integral curve of $\nabla f$. Since
\begin{equation}
    \frac{d}{ds}f(\gamma(s)) = df(\gamma'(s)) = df(\nabla f) = 1,
\end{equation}
it follows that $f(\gamma(s)) = f(\gamma(0)) + s$. Since $\gamma$ is a causal curve that is inextendible and $f$ is Cauchy temporal, it follows that $\lim_{s \uparrow b} f(\gamma(s)) = +\infty$, hence $b = +\infty$ as a consequence. Similarly, $a = -\infty$, thus $\gamma$ is complete. To conclude the proof, we consider the flow map $\mathrm{Fl}^{\nabla f}:\R \times \Sigma_0 \to M$. It is clearly a diffeomorphism as a consequence of what has been shown so far, and it is isometric as a consequence of $\mathrm{Hess} \, f = 0$.
\end{proof}

\begin{remark}[The case $p=1/2$]
Let us remark on some aspects of the proof of Theorem \ref{theorem: cosmological splitting p harmonicity} for $p=1/2$. In that case, the quantity $C_p(s)$ defined in \eqref{eq: def C_p} is precisely the hypersurface area
\begin{equation}
    \mathcal{A}(s):=\int_{\Sigma_s}d\mu_s
\end{equation}
of the $s$-level of $f$. Thus \eqref{eq: C_p' = p/q B_p} and \eqref{eq: derivative of B_p} give
\begin{equation}
    \mathcal{A}'(s)=B_{1/2}(s),
    \qquad
    \mathcal{A}''(s)
    =-\int_{\Sigma_s}\frac{Q}{w}\,d\mu_s\le0.
\end{equation}
In fact, the estimate \eqref{eq: B_p squared ineq} yields
\begin{equation}
    \mathcal{A}(s)\mathcal{A}''(s)
    +(\mathcal{A}'(s))^2\le0,
\end{equation}
so $\mathcal{A}^2$ is also concave. Since $\mathcal{A}$ is
positive on all of $\mathbb{R}$, it is constant. This is of course also clear a posteriori for any $p < 1$, once the splitting in Theorem \ref{theorem: cosmological splitting p harmonicity} has been established.
\end{remark}

\subsection*{Conflinct of interest statement}

The author confirms that there is no conflict of interest associated with this work.

\subsection*{Data availability statement}

The author confirms that there is no data associated with this work.

\subsection*{Use of AI declaration}

Discussions with ChatGPT 6 Astra were used to arrive at the precise proof strategies for the main results, and they were also helpful for calculations, exposition and proof-reading. The main ideas and final formulations of the manuscript are solely due to the author, and he takes full responsibility for the contents.

\subsection*{Acknowledgments}

The author would like to thank the Faculty of Mathematics of the University of Vienna for the hospitality during the conference ``Curvature, Optimal Transport and Gravity" in September 2026.

This research was funded in part by the Austrian Science Fund (FWF) [Grants DOI\\10.55776/EFP6, 10.55776/STA32 and 10.55776/J4913]. For open access purposes, the author has applied a CC BY public copyright license to any author accepted manuscript version arising from this submission.

\bibliography{references} 
\bibliographystyle{acm}

\end{document}